\documentclass[a4paper,11pt]{article}
\usepackage{latexsym}
\usepackage{amssymb}
\usepackage{amsfonts}
\usepackage{amsmath}
\usepackage{indentfirst}
\usepackage{graphicx}
\usepackage[noend,boxruled,vlined]{algorithm2e}
\usepackage{epsfig}

\newtheorem{prop}{Proposition}[section]
\newtheorem{teor}{Theorem}[section]
\newtheorem{lemma}{Lemma}[section]

\newcommand{\cvd}{\hfill $\blacksquare$\bigskip}

\date{}

\author{Luca Ferrari\thanks{Dipartimento di Matematica e Informatica
``U. Dini", viale Morgagni 65, 50134 Firenze, Italy {\tt
luca.ferrari@unifi.it}}}

\title{Greedy algorithms and poset matroids\footnote{Partially supported by PRIN
project: ``Automi e Linguaggi Formali: Aspetti Matematici e
Applicativi".}}

\begin{document}

\maketitle

\begin{abstract}
We generalize the matroid-theoretic approach to greedy algorithms
to the setting of poset matroids, in the sense of Barnabei,
Nicoletti and Pezzoli (1998) \cite{BNP}. We illustrate our result
by providing a generalization of Kruskal algorithm (which finds a
minimum spanning subtree of a weighted graph) to abstract
simplicial complexes.
\end{abstract}

\section{Introduction}

An \emph{independence system} is a pair $(E,\mathcal{F})$ such
that $E$ is a finite set and $\mathcal{F}$ is a down-set of the
Boolean algebra $\wp (E)$. A \emph{matroid} is an independence
system satisfying the following axiom: for any $A,B\in
\mathcal{F}$ such that $|B|=|A|+1$, there exists $b\in B\setminus
A$ such that $A\cup \{ b\} \in \mathcal{F}$.

In the paper \cite{BNP} the authors propose a generalization of
the notion of matroid where the ground set is equipped with a
partial order. The central definition of their work is the
following: a \emph{poset matroid} is a pair $(P,\mathcal{I})$
where $P$ is a finite partially ordered set and $\mathcal{I}$ is a
nonempty family of up-sets of $P$ satisfying the following
properties:
\begin{itemize}
\item[(i)] if $X,Y$ are up-sets of $P$ such that $Y\in
\mathcal{I}$ and $X\subseteq Y$, then $X\in \mathcal{I}$;

\item[(ii)] for every $X,Y\in \mathcal{I}$ with $|X|<|Y|$, there
exists $y\in Max(Y\setminus X)$ such that $X\cup \{ y\} \in
\mathcal{I}$.
\end{itemize}

The elements of $\mathcal{I}$ are called the \emph{independent
sets} of the poset matroid. To understand the above definition, we
recall that an up-set (resp., down-set) of a poset $P$ is a subset
$S$ of $P$ such that, if $x\in S$ and $x\leq y$ (resp. $x\geq y$),
then $y\in S$. Moreover, for any $S\subseteq P$, we denote by
$Max(S)$ the set of maximal elements of $S$.

We remark that the definition given here differs from the original
one in \cite{BNP}, which is given in terms of the notion of
\emph{basis}. However, the two definitions are clearly equivalent,
as it is shown in \cite{BNP}.

Given this generalized notion of matroid, it is natural to try to
generalize notions and results of matroid theory to the context of
poset matroids. Among the open problems proposed in \cite{BNP},
the last one is the following: is it possible to generalize the
generic greedy algorithm to the setting of poset matroids? To
better understand this problem, recall that there is a strong
relationship between greedy algorithms and the notion of matroid,
which we will briefly summarize below.

Given a \emph{weight function} $w:E\rightarrow \mathbf{R}^+$, we
consider the following problem:
\begin{align*}
\textnormal{\textbf{input: }}&\textnormal{an independence system
$(E,\mathcal{F})$ and a weight function}\\
&\textnormal{$w:E\rightarrow \mathbf{R}^+$.}\\
\textnormal{\textbf{output: }}&\textnormal{a set $M\in
\mathcal{F}$ such that $w(M)=\sum_{x\in M}w(x)$ is maximum.}
\end{align*}

The \emph{greedy algorithm} for the independence system
$(E,\mathcal{F})$ attempts to solve the above problem, and
consists of the following procedure:

\begin{algorithm}[H]
$S:=\emptyset$\; $Q:=E$\; \While{$Q\neq \emptyset$}{find $m\in Q$
having maximum weight\; $Q:=Q\setminus \{ m\}$\; \If{$S\cup \{ m\}
\in \mathcal{F}$}{$S:=S\cup \{ m\}$\;}} \Return{$S$}\;
\caption{GREEDY($(E,\mathcal{F}),w$) }\label{greedy}
\end{algorithm}

The procedure GREEDY tries to find a global solution by making the
local best choice at each step. Unfortunately, such an algorithm
is not always correct (that is, it does not solve the above
problem in general). The following theorem by Edmonds and Rado
\cite{E,R} tells us in which cases it works.

\begin{teor} Given an independence system $(E,\mathcal{F})$, the following
statements are equivalent:

\begin{itemize}
\item[a)] for any weight function $w$, GREEDY is correct on input
$(E,\mathcal{F})$, $w$;

\item[b)] $(E,\mathcal{F})$ is a matroid.
\end{itemize}
\end{teor}

In the next section we will generalize the Edmonds-Rado theorem to
the setting of poset matroid. In section  \ref{application} we
will see how our generalization can be used to find an analog of
Kruskal algorithm, which determines a minimum spanning subtree of
a weighted graph, where the graph is replaced by an abstract
simplicial complex. Finally, in the last section we will give some
hints to relate our work with previous approaches on the same (or
perhaps similar) subject.

\section{The Edmonds-Rado theorem for poset matroids}

Given a poset $P$, let $\mathcal{I}$ be a family of up-sets of $P$
satisfying condition (i) in the definition of poset matroid (i.e.
$\mathcal{I}$ is a down-set of up-sets of $P$). Call such a pair
$(P,\mathcal{I})$ a \emph{po-independence system}.

Consider the following problem:
\begin{align*}
\textnormal{\textbf{input: }}&\textnormal{a po-independence
system $(P,\mathcal{I})$ and a weight order-preserving}\\
&\textnormal{function $w:P\rightarrow \mathbf{R}^+$.}\\
\textnormal{\textbf{output: }}&\textnormal{an up-set $M\in
\mathcal{I}$ such that $w(M)=\sum_{x\in M}w(x)$ is maximum.}
\end{align*}

To solve it we try to adapt the greedy algorithm as follows:

\begin{algorithm}[H]
$S:=\emptyset$\; $Q:=P$\; \While{$Q\neq \emptyset$}{find a maximal
element $m\in Q$ having maximum weight\; $Q:=Q\setminus \{ m\}$\;
\If{$S\cup \{ m\} \in \mathcal{I}$}{$S:=S\cup \{ m\}$\;}}
\Return{$S$}\; \caption{PGREEDY($(P,\mathcal{I}),w$)
}\label{pgreedy}
\end{algorithm}

Our main result is the following generalization of the
Edmonds-Rado theorem for poset matroids.

\begin{teor}\label{main} Given a po-independence system $(P,\mathcal{I})$, the following
statements are equivalent:
\begin{itemize}
\item[a)] for any weight order-preserving function $w$, PGREEDY is
correct on input $(P,\mathcal{I})$, $w$;

\item[b)] $(P,\mathcal{I})$ is a poset matroid.
\end{itemize}
\end{teor}

\emph{Proof.}\quad $a)\Rightarrow b)$\quad Suppose that
$(P,\mathcal{I})$ is not a poset matroid. This means that there
exist $A,B\in \mathcal{I}$, with $|A|=k$ and $|B|=k+1$, such that,
for all $b\in Max(B\setminus A)$, $A\cup \{ b\}$ is an up-set but
$A\cup \{ b\} \notin \mathcal{I}$. Consider the weight function
$w:P\rightarrow \mathbf{R}^+$ defined as follows:
\begin{equation}
w(x)=\left\{ \begin{array}{lll} \alpha (>1)&\qquad x\in A
\\ 1&\qquad x\in B\setminus A
\\ 0&\qquad x\notin A\cup B
\end{array}\right. .
\end{equation}

We start by observing that $w$ is order preserving. Indeed, let
$x,y\in P$ such that $x\leq y$. If $x\in A$, then also $y\in A$
(since $A$ is an up-set), whence trivially $w(x)=w(y)=\alpha$. If
$x\in B\setminus A$, then clearly $y\in B$ (since $B$ is an
up-set); now, if $y\in A$, then $w(x)=1<\alpha =w(y)$, whereas, if
$y\notin A$, then $w(x)=1=w(y)$. Finally, if $x\notin A\cup B$,
then trivially $w(x)=0\leq w(y)$.

Now let $S\in \mathcal{I}$ be the solution provided by PGREEDY.
Depending on its cardinality, $S$ is a subset of $A$ or it
contains all elements of $A$ and some elements not in $B\setminus
A$. In fact, the elements of $A$ are the first ones that are
chosen by PGREEDY, since they have maximum weight (at each step,
PGREEDY will select a maximal element among the remaining ones in
$A$). In case all the elements of $A$ have already been selected,
it is possible that some (possibly all) of the elements of $P$ not
belonging to $B\setminus A$ are included in $S$. Denote by $C$ the
set of these elements ($C$ may also be empty). Observe that
PGREEDY cannot choose other elements, since, by hypothesis, $A\cup
\{ b\}\notin \mathcal{I}$, for all $b\in Max(B\setminus A)$ (and
so PGREEDY never enters $B\setminus A$). Now, set $t=|A\cap B|$,
we have
\begin{eqnarray*}
w(S)\leq w(A\cup C)=w(A)+w(C)=\alpha \cdot |A|=\alpha \cdot k \\
w(B)=w(B\setminus A)+w(A\cap B)=(k+1-t)+\alpha \cdot t.
\end{eqnarray*}

Thus, if we choose $1<\alpha <1+\frac{1}{k-t}$, we get
$w(S)<w(B)$, that is $S$ has not maximum weight, whence PGREEDY is
not correct in this case.

\medskip

$b)\Rightarrow a)$\quad Let $S=\{ b_1 ,b_2 ,\ldots ,b_n \}$ be the
solution provided by PGREEDY on input $(P,\mathcal{I})$, $w$, and
suppose that $w(b_1 )\geq w(b_2 )\geq \cdots \geq w(b_n )$. Now
consider $A=\{ a_1 ,a_2 ,\ldots ,a_m \} \in \mathcal{I}$, with
$w(a_1 )\geq w(a_2 )\geq \cdots \geq w(a_m )$.

We start by observing that $m\leq n$. Indeed, suppose $n<m$; then
(since $(P,\mathcal{I})$ is a poset matroid) there would exist
$a_j \in Max(A\setminus S)$ such that $S\cup \{ a_j \} \in
\mathcal{I}$. Moreover, for every up-set $R\subseteq S\cup \{ a_j
\}$, we would obviously have $R\in \mathcal{I}$, so $a_j$ should
belong to $S$, which is not.

Next we will prove that $w(a_ i)\leq w(b_i )$, for all $i=1,\ldots
,m$. Suppose it is not, and let $k$ be the minimum index for which
$w(a_k )>w(b_k )$. Notice that $D=\{ b_1 ,\ldots b_{k-1}\} \in
\mathcal{I}$, up to rearranging the elements of $S$. This can be
achieved without losing the property $w(b_1 )\geq w(b_2 )\geq
\cdots \geq w(b_n )$, since $w$ is order-preserving. The same
argument also shows that $\{ a_1 ,\ldots ,a_k \} \in \mathcal{I}$.
Now, since $|D|+1=|\{ a_1 ,\ldots ,a_k \} |$, there exists $a_j
\in Max(\{ a_1 ,\ldots ,a_k \} \setminus D)$ such that $D\cup \{
a_j \} \in \mathcal{I}$. But $w(b_k )\geq w(a_j )$ (since at the
$k$-th step PGREEDY chooses the element having maximum weight
among the remaining maximal ones) and $w(a_j )\geq w(a_k )$ (since
$j\leq k$), whence $w(b_k )\geq w(a_k )$, which is contrary to the
assumption.

The two above facts implies that $w(A)\leq w(S)$, and so that $S$
is indeed the correct solution, as desired.\cvd

\section{Acyclic subcomplexes of an abstract simplicial
complex}\label{application}

In order to illustrate our generalization of Edmonds-Rado theorem
to poset matroids, we propose a generalization of the well-known
Kruskal algorithm, which constructs a minimum spanning subtree of
a weighted graph.

\bigskip

Recall that an \emph{abstract simplicial complex} on a finite set
$X$ is a family $\mathcal{C}$ of subsets of $X$ such that, if
$F\in \mathcal{C}$ and $G\subseteq F$, then $G\in \mathcal{C}$
(i.e., a down-set of the powerset of $X$ partially ordered by
containment). Given $F\in \mathcal{C}$, we say that $F$ is a
\emph{face of dimension $i$} of $\mathcal{C}$ when $|F|=i$. The
set of all faces of dimension $i$ of $\mathcal{C}$ will be denoted
$\mathcal{F}_i$. Therefore, if the maximum dimension of a face of
$\mathcal{C}$ is $k$ (also called the \emph{dimension} of
$\mathcal{C}$), then $\mathcal{C}=\bigcup_{i=0}^{k}\mathcal{F}_i$.
Moreover, given $\mathcal{D}\subseteq \mathcal{C}$, we say that
$\mathcal{D}$ is a \emph{subcomplex} of $\mathcal{C}$ when it is
itself an abstract simplicial complex.

\bigskip

The faces of an abstract simplicial complex can be partially
ordered in a natural way by containment. However, to be consistent
with the theory we have developed in the previous sections, we
rather need to consider the dual order. Thus, given $F,G\in
\mathcal{C}$, we define $F\leq G$ whenever $G\subseteq F$. Observe
that a subcomplex of $\mathcal{C}$ is an up-set of
$(\mathcal{C},\leq)$.

\bigskip

Suppose that $\mathcal{C}$ is an abstract simplicial complex of
dimension $k$. For any $2\leq h\leq k$, we say that
$\mathcal{D}\subseteq \mathcal{C}$ is an \emph{$h$-cycle} of
$\mathcal{C}$ when:

\begin{itemize}

\item $\mathcal{D}\subseteq \mathcal{F}_h$;

\item for every $F\in \mathcal{D}$ and for every $x\in F$, there
exists precisely one face $H\in \mathcal{D}$ such that $F\cap
H=F\setminus \{ x\}$.

\end{itemize}

When an abstract simplicial complex does not have any $h$-cycles
it will be called \emph{$h$-acyclic}. Observe that an $h$-cycle of
a complex $\mathcal{C}$ is not a subcomplex of $\mathcal{C}$.
Moreover, given a face $F\in \mathcal{C}$ of dimension $h+1$, the
set of all faces of $F$ of dimension $h$ is an $h$-cycle, which
will be denoted $<F>$.

The following key lemma is central in the proof of our final
result.

\begin{lemma}\label{inters} Let $\mathcal{D}_1 ,\mathcal{D}_2$ be two $h$-cycles
of the abstract simplicial complex $\mathcal{C}$ such that
$\mathcal{D}_1 \cap \mathcal{D}_2 \neq \emptyset$. Then
$\mathcal{D}=(\mathcal{D}_1 \cup \mathcal{D}_2 )\setminus
(\mathcal{D}_1 \cap \mathcal{D}_2 )$ is an $h$-cycle of
$\mathcal{C}$ as well.
\end{lemma}

\emph{Proof.}\quad Obviously $\mathcal{D}\subseteq \mathcal{F}_h$.
Now take $F\in \mathcal{D}$ and $x\in F$; suppose moreover
(w.l.o.g.) that $F\in \mathcal{D}_1 \setminus
\mathcal{D}_2\subseteq \mathcal{D}_1$. Since $\mathcal{D}_1$ is an
$h$-cycle, there exists a unique $H\in \mathcal{D}_1$ such that
$F\cap H=F\setminus \{ x\}$. Moreover, it is clearly $F\setminus
\{ x\} =H\setminus \{ y\}$, for some $y\in H$. If $H\notin
\mathcal{D}_1 \cap \mathcal{D}_2$, then there is nothing else to
proof. Otherwise, if $H\in \mathcal{D}_1 \cap \mathcal{D}_2$, then
in particular $H\in \mathcal{D}_2$, whence there is a unique $G\in
\mathcal{D}_2$ such that $H\cap G=H\setminus \{ y\}$. Once again,
we also have that $H\setminus \{ y\} =G\setminus \{ z\}$, for some
$z\in G$. Observe that $G\notin \mathcal{D}_1$, since otherwise
there would exist two distinct faces in $\mathcal{D}_1$ whose
intersection with $H$ equals $H\setminus \{ y\}$ (namely $F$ and
$G$). Thus, in particular, $G\in \mathcal{D}$, and we have:
\begin{eqnarray*}
F\cap G&=&((F\setminus \{ x\} )\cup \{x \} )\cap G=((F\setminus \{
x\} )\cap G)\cup (\{ x\} \cap G)\\
& &(H\setminus \{ y\} )\cap G=H\cap G=H\setminus \{ y\}
=F\setminus \{ x\} .
\end{eqnarray*}

Finally, observe that $G$ is the unique face in $\mathcal{D}$
having the above property, since otherwise there would exist two
distinct faces in $\mathcal{D}_2$ whose intersection with $H$
equals $H\setminus \{ y\}$.\cvd

We will also use a result of \cite{BNP}, which asserts that
property (ii) in the definition of a poset matroid can be replaced
by a sort of ``local version". We report the precise statement in
the next lemma.

\begin{lemma} (\cite{BNP}) Let $\mathcal{I}$ a nonempty family of filters of a
poset $P$ satisfying property $(i)$ in the definition of poset
matroid. Then the following are equivalent:

\begin{itemize}

\item[(ii)] for every $X,Y\in \mathcal{I}$ with $|X|<|Y|$, there
exists $y\in Max(Y\setminus X)$ such that $X\cup \{ y\} \in
\mathcal{I}$;

\item[(ii')] for every $X,Y\in \mathcal{I}$ with $|Y|=1+|X|$ and
$|X|=1+|X\cap Y|$, there exists $y\in Max(Y\setminus X)$ such that
$X\cup \{ y\} \in \mathcal{I}$.

\end{itemize}

\end{lemma}

Given $2\leq h\leq k$, define $\mathfrak{I}_h =\{
\mathcal{D}\subseteq \mathcal{C}\; |\; \textnormal{$\mathcal{D}$
is a subcomplex of $\mathcal{C}$}$ $\textnormal{and does not
contain $h$-cycles}\}$.

\begin{prop} For any given $h$, $(\mathcal{C},\mathfrak{I}_h )$ is a poset matroid.
\end{prop}

\emph{Proof.}\quad First of all, it is clear that, if
$\mathcal{D}\in \mathfrak{I}_h$ and
$\widetilde{\mathcal{D}}\subseteq \mathcal{D}$ is a subcomplex,
then $\widetilde{\mathcal{D}}\in \mathfrak{I}_h$ as well.

To conclude the proof it will be enough to show that property
\emph{(ii')} of the above lemma holds. So let $\mathcal{D}_1
,\mathcal{D}_2 \in \mathfrak{I}_h$ such that
$|\mathcal{D}_2|=1+|\mathcal{D}_1 |$ and $|\mathcal{D}_1
|=1+|\mathcal{D}_1 \cap \mathcal{D}_2 |$. Observe that, in this
situation, it is $|\mathcal{D}_1 \setminus \mathcal{D}_2 |=1$ and
$|\mathcal{D}_2 \setminus \mathcal{D}_1 |=2$. Suppose that
$\mathcal{D}_1 \setminus \mathcal{D}_2 =\{ x\}$ and $\mathcal{D}_2
\setminus \mathcal{D}_1 =\{ y_1 ,y_2 \}$. There are of course two
distinct possibilities concerning $y_1$ and $y_2$. Suppose first
that $y_1$ and $y_2$ are incomparable. By way of contradiction,
suppose there exist $h$-cycles $\mathcal{Z}_1 ,\mathcal{Z}_2$ such
that $\mathcal{Z}_i \subseteq \mathcal{D}_1 \cup \{ y_i \}$, for
$i=1,2$. Since $\mathcal{D}_2 \in \mathfrak{I}_h$, there are $x_1
,x_2 \in \mathcal{D}_1 \setminus \mathcal{D}_2$ such that $x_i \in
\mathcal{Z}_i$, for $i=1,2$. However, our hypotheses imply that
$x_1 =x_2 =x$, whence we would have (from lemma \ref{inters}) that
$(\mathcal{Z}_1 \cup \mathcal{Z}_2 )\setminus \{ x_1 \} \subseteq
\mathcal{D}_2$ contains an $h$-cycle, which is forbidden since
$\mathcal{D}_2 \in \mathfrak{I}_k$. Finally, suppose that $y_1
<y_2$. Once again, we argue by contradiction, supposing that there
exists an $h$-cycle $\mathcal{Z}\subseteq \mathcal{D}_1 \cup \{
y_2 \}$. This implies that $y_2$ is a face of dimension $h$, and
so $y_1\in \mathcal{D}_2$ has dimension $h+1$. Observe that all
the faces of $y_1$ but $y_2$ must be both in $\mathcal{D}_2$
(since $\mathcal{D}_2$ is an up-set) and in $\mathcal{D}_1$ (since
$\mathcal{D}_2 \setminus \mathcal{D}_1 =\{ y_1 ,y_2 \}$), whence
$(\mathcal{Z} \cup <y_1 >)\setminus \{ y_2 \} \subseteq
\mathcal{D}_1$. Moreover $\mathcal{Z}$ and $<y_1
>$ are $h$-cycles both containing $y_2$, hence, by lemma
\ref{inters}, $(\mathcal{Z} \cup <y_1 >)\setminus \{ y_2 \}$
contains an $h$-cycle, which is impossible.\cvd

We are now in a position to provide a Kruskal-like algorithm to
find a maximum spanning subcomplex of an abstract simplicial
complex with respect to a suitable weight function of its faces. A
\emph{spanning subcomplex} of a complex $\mathcal{C}$ is a
subcomplex of $\mathcal{C}$ containing all its 0-dimensional
faces. The weight of a (sub)complex is simply the sum of the
weights of its faces.

\begin{teor} Let $\mathcal{C}$ be an abstract simplicial
complex, and let $w: \mathcal{C}\rightarrow \mathbf{R}^+$ be an
order-reversing function (given that $\mathcal{C}$ is partially
ordered by containment). Then the algorithm PGREEDY is correct on
input $((\mathcal{C},\mathfrak{I}_h ),w)$, and returns a spanning
$h$-acyclic subcomplex of $\mathcal{C}$ having maximum weight.
\end{teor}

\section{Conclusions}

In this note we have extended the classical Edmonds-Rado theorem
to the more general setting of poset matroids described in
\cite{BNP}. We have illustrated our result by generalizing a
classical algorithm on graphs due to Kruskal to the setting of
abstract simplicial complexes. Of course, lots of other possible
applications can be considered. One of the most interesting is
perhaps the generalization of the greedy solution of the classical
task scheduling problem presented, for instance, in \cite{CLRS}.
The obvious modification of this very well-known application of
Edmonds-Rado theorem consists of introducing a \emph{priority}
between tasks, which can be naturally formalized as a partial
order relation. However, our attempts to find a correct analog of
this problem (and its solution) in the context of poset matroids
have been unsuccessful, so it would be very interesting to have
some results in this direction.

We remark that the extension of the concept of matroid on finite
sets to posets considered in the present paper is not the only one
that can be found in the literature. Another well known approach
is through the theory of \emph{geometries on partially ordered
sets} due to Faigle \cite{F2}, which is however intimately related
to the one proposed in \cite{BNP}.

Even more interestingly, in \cite{F1} Faigle finds a necessary and
sufficient condition for a generic greedy algorithm to be correct
in a setting that is extremely similar to ours. Apart from the
fact that he considers independent set to be down-set rather than
up-sets, which is an immaterial difference (it just consists of
dualizing all the definitions given here), the analogies with our
results are really striking. However, the conditions found by
Faigle (which are condensed in what he calls a ``generating set")
are slightly different from our, and it is not immediately evident
how to relate the two approaches.

Another well-known generalization of matroid theory, which is more
oriented towards greedy algorithms, is the theory of greedoids
introduced by Korte and Lov\'asz in \cite{KL}. In \cite{LZ} the
authors try to merge the notions of poset matroid and of greedoid
by developing the theory of poset greedoids. It would be
interesting to have a generalization of our results to the setting
of poset greedoids.

We conclude by recalling that in \cite{S} the author proves that
the correctness of a general greedy algorithm for a hereditary
system is equivalent to the fact that such system is a so-called
\emph{strict cg-matroid}. It is likely that there is a
relationship between the results of the present paper and those of
\cite{S}, but it is not clear to us how to make it explicit.

\end{document}